
\input amstex
\def\cal{\Cal}
\magnification=\magstep1
\documentstyle{amsppt}
\topmatter
\title NODAL CURVES ON SURFACES OF GENERAL TYPE \endtitle
\author L.Chiantini\\  E.Sernesi\endauthor
 \affil
    Dipartimento Me.Mo.Mat, Universita' 'La Sapienza', \\  via Scarpa 10, 00161
Roma (Italy) \\
e-mail: chiantin\@itcaspur.caspur.it \\ \\
Dipartimento di Matematica, Terza Universita' di Roma,  \\ via C.Segre 2, 00146
Roma (Italy) \\
          e-mail: sernesi\@matrm3.mat.uniroma3.it \endaffil

\endtopmatter

\document

{\bf INTRODUCTION}

In this paper we investigate to which extent
 the theory of Severi on nodal plane curves of a given degree $d$
extends to a  linear system on a complex projective
nonsingular algebraic surface. As well known, in  [S], Anhang
F Severi proved that for every $d \ge 3$ and
$0 \le \delta \le {d-1 \choose 2}$ the family $\cal V_{d,\delta}$ of plane
irreducible curves of degree $d$ having exactly $\delta$ nodes and no
other singularities is non empty and everywhere smooth of codimension
$\delta$ in the linear system $|{\cal O}(d)|$.

If $C \in \cal V_{d,\delta}$  Severi uses the non speciality of
the normal line bundle to the composition $\nu: \tilde C \to C\to \bold P^2$,
where $\tilde C$
is the normalization of $C$, to prove that $\cal V_{d,\delta}$ is smooth of the
asserted codimension at the point $ C$. This proof can  be
extended to rational, ruled and  K3 surfaces with little changes: we
discuss this point in section 1 (see also [T] for the case of rational
surfaces.

In the case of a surface of general type $S$ the  approach of Severi fails,
and in fact it is easy to see that the analogous of Severi's theorem does
not hold in general if we  impose too many
nodes to the curves of a complete linear system $|D|$.  One may
nevertheless look for  an upper bound on $\delta$ ensuring that the family
 $\cal V_{D,\delta}$   of irreducible curves in $|D|$ with $\delta$ nodes is
smooth of codimension $\delta$. In section 2 we  give the following partial
answer to
this problem:

\proclaim{ Theorem 2.2}    Let $S$ be a surface such that $|K_S|$ is ample
and let $C$ be an irreducible curve on $S$ such that
$C \equiv _{num} pK_S$, $p\ge 2$, $p\in {\bold Q}$ and $|C|$ has smooth general
member.
Assume that $C$ has
$\delta \ge 1$ nodes and no other singularities and assume that either
$ \delta < {p( p-2) \over 4}  K^2_S, $
or $ \delta < {(p-1)^2 \over 4}  K^2_S $, $p\in {\bold Z}$ odd, and the
Neron-Severi group of $S$ is ${\bold Z} $ generated by $K_S$.

Then $\cal V_{C,\delta}$ is smooth of  codimension $\delta$
at the point corresponding to $  C $.\endproclaim

For the proof  we consider a curve
 $  C   \in \cal V_{D,\delta}$ where  Severi's theorem fails
and we associate
a rank two vector bundle  to the zero-cycle $N$ of
nodes of $C$. We then apply
 the Bogomolov-Reider
method ([Rr]) to deduce the inequality of the theorem from the properties of
this
vector bundle.

As a special case of theorem 2.2 we obtain the following
result on smooth
 surfaces in ${\bold P}^3$:

\proclaim{Proposition 2.4}   Let $S$ be a smooth surface of degree $d\ge 5$ in
 ${\bold P}^3$ with plane section $H$. If $C\in |nH|$ has $\delta$ nodes and
no other singularities and $\delta < {nd(n-2d+8)\over 4}$ then
 $\cal V_{C,\delta}$ is smooth of  codimension $\delta$
at the point $  C $.\endproclaim

In the case of quintic surfaces we slightly improve this bound
proving that $\cal V_{C,\delta}$ is smooth of  codimension $\delta$
at the point $  C $ if

 $\delta < {5(n-1)^2 \over 4}$
 when $n$ is odd (proposition 2.5).
 We show that  these estimates are  sharp,
by producing explicit examples of  curves in the linear
systems $|nH|$ having exactly $\delta$ nodes, where $\delta$ is the upper
bound given above and which are obstructed as elements of $\cal V_{C,\delta}$.
These examples are  discussed in  detail in section 4,
 where we show in
particular that they are not general points of a component of $\cal
V_{C,\delta}$.

The construction
of the
 examples is achieved through the consideration of another
problem. We consider a complete intersection curve $C$ in ${\bold P}^3$
having $\delta$ nodes;  call $C$  {\sl geometrically linearly
normal} if it is not a birational projection of a smooth curve of
 ${\bold P}^4$ of the same degree.
Since clearly smooth complete intersections are geometrically linearly normal,
 one may look for a  bound $\delta (d,n)$ such that $C$ is geometrically
linearly
normal if $\delta \le \delta (d,n)$.  Aiming
at this we  prove the
following:
\proclaim{ Theorem 3.4}  Assume that $C\subset\bold P^3$ is the complete
intersection of
a smooth surface $S$ of degree $d$ with a surface of degree $n$.  $C$ is
geometrically linearly normal if and only if the nodes of $C$ impose
independent conditions to the linear system  $|(n+d-5)H|$ on $S$.\endproclaim

As a conseguence of this criterion we find the upper bound
$\delta (d,n) = {nd(n-2) \over 4}$  (theorem 3.5).
Theorem 3.4 is applied to the case of curves
on a general quintic surface to
deduce that certain non geometrically linearly
normal curves  we construct  are obstructed in the corresponding Severi
variety.
  These examples also show that the previous bound   $\delta (5,n)$ is sharp.

We have not considered the problem of existence for  $\cal V_{D,\delta}$;
 for results in this direction we refer to [CR] and [X].

The paper consists of four section. Section 1 is devoted to
known facts  on nodal curves and to Severi theory on  rational, ruled
 and
K3 surfaces. In section 2 we prove our main results on surfaces of general
type. Section 3 deals with geometric linear normality of nodal complete
intersection curves in  ${\bold P}^3$. In section 4 we construct
the examples of nodal curves  on a quintic surface which show that the
results of section 2  are sharp.

We  work in the category of schemes over ${\bold C}$, the field of
complex numbers. As usual, $\dim(H^i($--$))$  will be denoted by
$h^i($--$)$.

\bigskip

{\bf 1. PRELIMINARIES}
\medskip

We will denote by $S$ a  projective nonsingular algebraic
surface. Let $|D|$ be a complete linear system on $S$ whose general member is
an irreducible
non singular curve.
 We will denote by $p_a(D)$ the {\it arithmetic genus}
of $D$, given by:
$$p_a(D) = {D(D+K_S) \over 2}+1.$$
For every $\delta \ge 0$ there is a locally closed
subscheme $\cal{V}_{D,\delta}$ of $|D|$ which parametrizes a
universal family of reduced and irreducible curves belonging to
$|D|$ and having exactly
$\delta$ nodes (ordinary double points)
and no other singularities (see [W] for $S=\bold P^2$, but the proof extends to
any $S$).\par
The schemes $\cal{V}_{D,\delta}$ will be called {\it Severi varieties}.
Let $ C \in{\cal V}_{D,\delta}$ and let $N$ be the scheme of nodes
of $C$: it is a closed zero-dimensional subscheme of $S$ of
degree $\delta$. The {\it geometric genus} of $C$ is
$$g = p_a(D) - \delta$$
 The Zariski tangent space of $|D|$ at $C$ is
$$T_{|D|,C}={H^0(S,{\cal O}_S(D))\over (C)}$$
and the Zariski tangent space of $\cal{V}_{D,\delta}$  at $C$ is
$$T_{{\cal V}_{D,\delta},C}= {H^0(S,{\cal I}_N(D))\over (C)}$$
while the obstruction space is a subspace of $H^1(S,{\cal I}_N(D))$.\par
In other words, a first order deformation
$C +\epsilon C'$, $\epsilon^2 = 0$, is in  ${\cal V}_{D,\delta}$ if and
only if it is in $|D|$ and $N\subset C'$. In particular:
$$ \dim(T_{{\cal V}_{D,\delta},C })\ge h^0(S,{\cal O}_S(D)) - \delta - 1
 = h^0(S,{\cal O}_S(D)) - (p_a(D) - g) - 1 $$
and equality
holds iff $N$ imposes independent conditions to $|D|$. In
this case ${\cal V}_{D,\delta}$ is nonsingular of dimension
$$h^0(S,{\cal O}_S(D)) - \delta - 1 = \dim(|D|) - \delta$$
at $C  $.

We recall the theorem of Severi on the projective plane.

\proclaim {Theorem 1.1} ({\rm Severi}) Let $S = {\bold P}^2$,
$d\ge 3$ and $D$ any divisor of degree $d$. Let $\delta
\ge 1$ be such that  $$ \delta \le  p_a(D) = {d-1\choose 2}.$$
Then the Severi variety ${\cal V}_{D,\delta}$ is non empty and smooth of
pure
 dimension  $$\dim(|D|) - \delta = {d(d+3)\over 2} - \delta.$$
\endproclaim
\demo{Proof}  Let us suppose that $C\in {\cal V}_{D,\delta}$ and let $N$
be the scheme of nodes of $C$. In view of the exact sequence
$$0\to {\cal I}_N(d)\to {\cal O}_S(d)\to {\cal O}_N(d)\to 0  $$
and of the fact that $H^1(S,{\cal O}_S(d))=0$, in order to prove that
${\cal V}_{D,\delta}$  is smooth of the asserted dimension at $C$ it
is necessary and sufficient to prove that $H^1(S,{\cal I}_N(d))=0$.\par
Let $\sigma :=h^
1(S,{\cal I}_N(d))$. Since  $h^0({\cal O}_N(d))=\delta$,
from the above sequence we deduce that
$$h^0(S,{\cal I}_N(d)) = h^0({\cal O}_S(d))-\delta+\sigma={d+2\choose
2}-\delta+\sigma.$$
Let $\nu :\tilde C \to C$ be the normalization of $C$ and let
$\tilde N$ be the pullback of $N$ to $\tilde C$. We have an injective map:
$${H^0(S,{\cal I}_N(d))\over (C)}\to H^0(\tilde C,\nu ^*\cal O(d)(-\tilde N))
\tag 1$$
Since $\nu ^*\cal O(d)(-\tilde N))$ has degree
$$d^2-2\delta=2g-2+3d$$
then it is non-special
 and we deduce that:
$$  h^0(S,{\cal I}_N(d))-1\le h^0(\tilde C,\nu ^*\cal O(d)(-\tilde N))
=d^2-2\delta +1-g={d+2\choose 2}-\delta -1, $$
whence $\sigma =0$.\par
To prove that ${\cal V}_{D,\delta}\ne \emptyset$ for all $\delta$ we
start from the case $\delta =p_a(D)$ i.e. $g=0$.\par
The family ${\cal V}_{D,p_a(D)}$ is not empty because it contains any
general projection of a rational and normal curve of ${\bold P}^d$. Let
$C\in {\cal V}_{D,p_a(D)}$, let as usual $N$ denote the scheme of
nodes of
$C$, let $P\in N$ and $M$ the complement of $P$ in $N$. Since
$$h^1(S,{\cal I}_N(d)) = h^1(S,{\cal I}_M(d))$$
we have
$$h^0(S,{\cal I}_M(d))=h^0(S,{\cal I}_N(d))+1.$$
Any element of the vector space $H^0(S,{\cal I}_M(d))$ not in
$H^0(S,{\cal I}_N(d))$ defines an  infinitesimal
deformation of $C$ which smooths the node $P$ while leaving unsmoothed
all the other nodes. This means that $C\in \overline{{\cal V}}_{D,p_a(D)-1}$,
the closure of ${\cal V}_{D,p_a(D)-1}$.  Therefore ${\cal V}_{D,p_a(D)-1}\ne
\emptyset$. By descending induction on $\delta$ one proves similarly that
$${\cal V}_{D,\delta}\ne \emptyset$$
for all $1\le \delta \le p_a(D)$. \qed
\enddemo

In  example  1.3  we show how the proof of theorem  1.1  can be
adapted to K3 surfaces.

\proclaim{Remark 1.2} {\rm The reason why the proof of 1.1 works is because
$$\nu^*\cal O(d)(-\tilde N)=\nu^*(\cal O(d-3)(-\tilde N) \otimes\cal O(3))=
K_{\tilde C}\otimes\nu^*\cal O(3)$$
and therefore this line bundle is non special. This fact has
been applied in
\thetag2 to get $\sigma=0$.\par
It is then clear that if we consider any rational or ruled surface $S$ and any
smooth and
irreducible curve $C$ on $S$, such that $|C|$ is base point free and $K_SC<0$,
then
the first part of the proof of 1.1 (excluding the existence statement) can be
repeated to
this case word by word; this holds in particular for any Del Pezzo surface. For
this
result we refer also to [T]. Therefore we get the  following:\par\smallskip
\sl Let $S$ be a rational or ruled
 surface and let $C\subset S$ be a smooth irreducible curve
such that $|C|$ is base point free and $K_SC<0$. If for some $\delta\le p_a(C)$
we have
$\cal V_{C,\delta}\neq\emptyset$, then $\cal V_{C,\delta}$ is smooth of
codimension $\delta$
in $|C|$.\rm\endproclaim

\proclaim{Example 1.3} \rm Let $S$ be a K3 surface and $D$ a smooth irreducible
curve
such that $p_a(D)\ge 2$. Then (see [M]) $|D|$ is base point free and of
dimension $p_a(D)$; moreover $H^1(S,{\cal O}_S(D))=0$. For
each  $1\le
\delta \le p_a(D)$ and for any $C\in {\cal V}_{D,\delta}$ we have
(with notations as above):
$$\multline h^0(S,{\cal I}_N(C))-1 = p_a(D)-\delta +h^1(S,{\cal I}_N(C))\le\\
\le h^0(\tilde C,\nu ^*\cal O(C)(-\tilde N))=h^0(\tilde C,K_{\tilde
C})=p_a(D)-\delta.
\endmultline $$
It follows that $H^1(S,{\cal I}_N(C))=0$ and therefore
${\cal V}_{D,\delta}$ is smooth and of codimension $\delta$ in $|D|$.\par
In [MM] it is shown that ${\cal V}_{D,p_a(D)}\ne \emptyset$: therefore
as in the proof of 1.1 it follows
 that
${\cal V}_{D,\delta}\ne \emptyset$ for all $1\le \delta \le p_a(D)$.\par
Note that in particular we have that ${\cal V}_{D,p_a(D)}$ is finite,
i.e. {\it there are finitely many nodal rational curves in} $|D|$.
\endproclaim

\bigskip

{\bf 2. SURFACES OF GENERAL TYPE}
\medskip

If $S$ is a surface of general type then we cannot expect
that theorem  1.1  extends without changes to linear systems on $S$.
The reason for this is obvious. If $|D|$ is a (say very ample) linear
system on $S
$, then on a general curve $C\in |D|$ the characteristic
linear series is special; this implies that
$$\dim(|D|)\le g(C)-1=p_a(D)-1$$
therefore ${\cal V}_{D,p_a(D)}$ cannot have the expected codimension
and we should in fact expect that  ${\cal V}_{D,p_a(D)}=\emptyset$.

In this case we should ask the following more appropriate:

\proclaim{(2.1) Question } Given a surface of general type $S$ and a linear
system $|D|$ on $S$ whose general member is smooth and connected, for
which values of $\delta$
 is ${\cal V}_{D,\delta}$ non empty and smooth
of codimension $\delta$?
\endproclaim

We will give a partial answer to question (2.1). Our main result is the
following:

\proclaim{Theorem 2.2} Let $S$ be a surface such that $|K_S|$ is ample, and
let $C$ be an irreducible curve on $S$ such that $|C|$ contains smooth elements
and such that
$$C\equiv _{num}pK_S \qquad p\ge2,\  p\in {\bold Q}$$
Assume that $C$ has $\delta \ge 1$ nodes and no other singularities and assume
that either
 $$\delta < {p(
p-2) \over 4} K_S^2 $$
or
 $$\gather\delta < {(p-1)^2 \over 4} K_S^2\quad p\in\bold Z\text{ odd,   and
the Neron Severi}\\
\text { group of $S$ is  generated by } K_S.\endgather$$
Then the nodes of $C$ impose independent conditions to $|C|$. In
particular the Severi variety ${\cal V}_{C,\delta}$ is smooth of
codimension $\delta$ at  C.  \endproclaim

The proof is based on the study of a rank 2 bundle on $S$, associated with the
set
of nodes of $C$. To do this, we recall briefly the connections between
 rank 2 bundles
and sets of points on a surface.

\proclaim{Remark 2.3} \rm (see  [GH]) Let $N$ be a set of $\delta$ points in
$S$. If $N$ does not impose independent conditions to a linear system $|C|$,
then the restriction map $H^0(S,\cal O_S(C))\to H^0(\cal O_N)$ is not
surjective. Let $N_0\subset N$ be a minimal subset for which the composition
$H^0(S,\cal O_S(C))\to H^0(\cal O_N)\to H^0(\cal O_{N_0})$ does not surject;
then a general element of $H^1(S,\cal I_{N_0}(C))$ defines an extension:
$$
 0\to \cal O_S\to E \to \cal I_{N_0}(C-K_S)\to 0   \tag 2 $$
where $E$ is a rank 2 vector bundle on $S$, with Chern classes
$$\gather  c_1(E)=\cal O_S(C-K_S) \\
           c_2(E)= \deg N_0 \le\delta   \endgather $$
with also $c_2(E)>0$ for $N_0$ cannot be empty.
\endproclaim

\demo{Proof of Theorem 2.2} Call $N$ the set of nodes of the curve $C$ and
assume that
$N$ does not impose independent conditions to the curves of $|C|$; we show
that we get a contradiction. \par
Take the subset $N_0\subset
 N$ and the rank 2 vector bundle $E$ described in the previous remark and
denote by $\delta_0 $ the degree of $N_0$.\par
By assumptions, $c_1(E)=\cal O_S(C-K_S)\equiv_{num} (p-1)K_S$ and
$$ c_1(E)^2-4c_2(E) = (p-1)^2K_S^2-4\delta_0\ge (p-1)^2K_S^2-4\delta>0 $$
so that $E$ is {\it Bogomolov unstable} (see [B]).\par
It follows that there exists a divisor $M$ which 'destabilizes'  $E$ with
respect to the ample
divisor $K_S$, that is, $ h^0(S,E(-M))>0 $ and
$$(2M-c_1(E))K_S >0 \quad\text{i.e. }MK_S>({p-
1\over 2})K_S^2. \tag 3 $$
Taking $M$ maximal, we may further assume that a general section of $E(-M)$
vanishes in a locus $Z$ of codimension 2 (see [R] th.1).
It follows $\deg Z=c_2(E(-M))\ge 0;$ hence:
$$ \delta_0+M^2-(p-1)MK_S = c_2(E)+M^2-Mc_1(E) = c_2(E (-M)) \ge 0. \tag 4 $$
Let us now use \thetag 2. $h^0(S,\cal O_S(-M))$ is $0$,
for $-MK_S<-(p-1)K_S^2/2\le 0$ by assumptions and $K_S$ is ample; thus
$h^0(S,E(-M))>0$ implies $h^0(S,\cal I_{N_0}(C-K_S-M)>0$, that is, there exists
a
divisor $
\Delta$ in the linear system $|C-K_S-M|$, which contains  $\delta_0$ nodes
of the curve $C$. $\Delta$ cannot contain $C$ as a component, for as above
$(-K_S-M)K_S<0$, hence $-K_S-M$ cannot be effective. It follows, by Bezout,
$(C-K_S-M)C\ge 2\delta_0$ which yields:
$$ ((p-1)K_S-M)(pK_S)\ge 2\delta_0. \tag 5$$
Now observe that, since $K_S$ is ample, by Hodge index theorem, we have
$M^2K_S^2\le (MK_S)^2$;
putting this together with \thetag 4 and \thetag 5, one finally gets:
$$ \frac{(MK_S)^2}{K_S^2}-
{3p-2\over 2}MK_S+{p(p-1)\over 2}K_S^2\ge 0 $$
 which can be solved with respect to $MK_S$. Since by assumption $p\ge 2$  we
have
$p-1\ge p/2$, thus  the previous inequality implies either $MK_S\le pK_S^2/2$
or $MK_S\ge(p-1)K_S^2$. The last inequality yields
$((p-1)K_S-M)K_S\le 0$ while $(p-1)K_S-M\equiv_{num} C-K_S-M$ and $|C-K_S-M|$
has the
effective divisor $\Delta$ which contains $N_0\neq\emptyset$: this is
impossible
for $K_S$ is ample.\par
It remains to exclude $MK_S\le pK_S^2/2$. If the Neron
 Severi group of $S$ is $\bold Z$, generated by $K_S$ and $p$ is odd, then the
intersection of any two divisors on $S$ is an integral multiple of $K_S^2$, so
this inequality  implies  $MK_S\le (p-1)K^2_S/2$, which is  excluded by
\thetag3. Otherwise, use the assumption $\delta_0\le\delta<p(p-2)K_S^2/4$ in
\thetag4, together with Hodge inequality; we get:
$$ \frac{(MK_S)^2}{K_S^2}-(p-1)MK_S+{p (p-2) \over 4}K_S^2>0 \tag6$$
from which it follows that either $MK_S<(p-2)K^2_S/2$, absurd by \thetag 3, or
$
MK_S>{p\over 2}K_S^2  $, which yields the required contradiction.
\qed\enddemo

  We will show that our estimate on $\delta$ for having ${\cal V}_{C,\delta}$
smooth, of the expected codimension, is in fact sharp at least in some
example.\par
Let us point out that, for surfaces in $\bold P^3$ of degree $d\ge 5$, we have
$K_S=(d-4)H$  (very) ample and we may apply
the theorem to any curve $C\in |nH|$, for $n\ge 2(d-4)$, getting:

\proclaim{Proposition 2.4}   Let $S$ be
a smooth surface of degree $d\ge 5$ in $\bold P^3$ with
plane section $H$. If $C\in |nH|$, $n\ge 2d-8$ has $\delta$ nodes and no other
singularities, and
$\delta < nd(n-2d+8)/4$, then $C$ corresponds to a smooth point of a component
of the
Severi variety ${\cal V}_{C,\delta}$ with the expected codimension $\delta$.
\endproclaim

When $S$ is a general   quintic surface in $\bold P^3$ and $C=pK_S=pH$, $p$ odd
integer, theorem 2.2  gives:

\proclaim{Proposition 2.5}   Let $S$ be a   smooth surface of
 degree $ 5$ in $\bold P^3$
with plane section $H$ and Picard group $\bold Z$.
If $C\in |pH|$ ($p\ge 3$ and odd)   has $\delta$ nodes and no other
singularities, and
$\delta < 5(p-1)^2/4$, then  the Severi variety ${\cal V}_{C,\delta}$  is
smooth, with the
(expected) codimension $\delta$ at $C$.
\endproclaim

\proclaim{Remark 2.6} \rm One may apply the previous procedure also to K3 or
rational surfaces and get  estimates on  $\delta$ which implies that  ${\cal
V}_{C,\delta}$
is smooth, of the exp
ected codimension. However, for these surfaces
we get statements which are weaker than theorem 1.1 or example 1.2.\par
On the other hand, when $S$ is any smooth 5-ic surface in $\bold P^3$, we are
going to
provide examples that show that the numerical bounds for $\delta$ found in
theorem 2.2
and proposition 2.5  are sharp.\endproclaim

\proclaim{Remark 2.7} \rm Let $S$, $C$, $p$, $\delta$ be as in the statement of
theorem 2.2, but assume now:
$$\delta = {p (p-2 )\over 4}K_S^2.$$
If the nodes of $C$
 do not impose independent conditions to $|C|$, then we may
go through the proof of the theorem, finding the rank 2 bundle $E$ associated
to a subset
$N_0\subset N$ and the destabilizing divisor $M$. The only difference is that,
in
\thetag 3 , one gets only the weak inequality $MK_S\ge pK_S^2 /2$
so that, at the end of the argument, we cannot exclude the case
$MK_S=pK_S^2/2$.\par
If the equality holds, one deduces from \thetag 4 and \thetag 6:
$$ c_2(E(-M))= \delta_0-{p (p-2 )\over 4}K_S^2=\delta_0-
\delta.$$
  Since $E(-M)$ has sections vanishing in codimension 2, then   $c_2(E(-M))=0$,
so  $N_0=N$ and $E(-M)$ must split, thus
$$  E= \cal O_S(M)\oplus\cal O_S(C-M)$$
and $N$ is complete intersection of type $M, C-M$ on $S$.
\endproclaim
\bigskip

{\bf 3. A  RELATED PROBLEM: GEOMETRICALLY LINEAR NORMALITY}
\medskip

Let $C$ be a smooth, complete intersection curve in $\bold P^3$. It is well
known
that $C$ is arithmetically normal, so that, in particular, $C$ cannot be the
birational
proje
ction of a non-degenerate curve $C'\in\bold P^r$ for $r>3$, that is, the
embedding
$C \to \bold P^3$ does not factorize through any non-degenerate map $C \to
P^r$, $r>3$.\par

When $C$ has singularities, this is no longer true (as we shall see later):
there are
complete intersection singular curves $C\in\bold P^3$ whose normalization
$\tilde C\to C$ factors through a birational non degenerate map $\tilde C\to
\bold P^r$ for some
$r>3$. On the other hand, when the geometric genus of $C$ is close enough
 to the arithmetic
genus  of $C$, this factorization is impossible.  So, one may look for bounds,
for
the number $p_a(C)-g$, which exclude that $C$ can be obtained as the birational
projection of
a non-degenerate curve lying in some higher dimensional projective space.\par
In fact, we shall look at the case of curves having only nodes for
singularities
and lying on some fixed smooth surface $S$.

\proclaim{Definition 3.1} \rm Let  $C$ be any reduced curve in $\bold P^r$. We
say that $C$ is {\it '
geometrically linearly normal'} if the normalization $\tilde C\to C$
cannot be factored with a birational non-degenerate map $\tilde C\to \bold
P^R$, $R>r$,
followed by a projection.\endproclaim

\proclaim{(3.2) Problem} Let $S$ be a smooth surface of degree $d$ in $\bold
P^3$;
for any number $n$ find a sharp bound $\delta(d,n)$ such that if
$C\subset S$ is a complete intersection curve of type $d,n$, having only
$\delta$ nodes
as singularities and $\delta\le\delta(d,n)$, then $C$ is geometrically
linearly normal.
\endproclaim

A partial answer to this problem can be still given using Reider's construction
as in the
proof of theorem 2.2 and it turns out that, on a quintic surface, question
(2.1) and
problem (3.2) are in fact closely related.\par
To begin with, let us recall the following, well-known fact:

\proclaim{Proposition 3.3}    Let $S$ be a smooth  surface of and let $H$ be a
very ample
 divisor on $S$, such that for all $m$ $h^1\cal O_S(mH)=0$; let $C\in|nH|$ be
an irreducible
 curve,
 having only nodes for singularities; call $N$ the set of nodes of $C$,
$\nu:\tilde C\to C$ the
normalization and $\tilde N$ the pull-back of $N$ on $\tilde C$.\par
For all integers $m$ we have an isomorphism:
$$ \frac{H^0(S,\cal I_N(mH+K_S))}{H^0(S,\cal I_C(mH+K_S))}\to
H^0(\tilde C,\nu^*\cal O(mH+K_S)(-\tilde N))$$\endproclaim
\demo{Proof} Call $\mu:\tilde S\to S$ the blowing up of $S$ along $N$ and let
$B=\sum E_i$ be the
exceptional divisor; then $\tilde C$ is
isomorphic to (and will
be identified with) a divisor on  $\tilde S$ in the class
$\mu^*C-2B$. Since $\omega_{\tilde C}$ is cut by the divisors in
$|\mu^*C+\mu^*K_S-B|$,
the exact sequence:
$$\multline 0\to\cal O_{\tilde S}(\mu^*K_S+B+(m-n)\mu^*H)\to\cal O_{\tilde
S}(\mu^*C+\mu^*K_S-B+
(m-n)\mu^*H)\to \\ \to\omega_{\tilde C}((m-n)\mu^*H)\to 0\endmultline$$
shows that the statement follows once we know that $h^1 O_{\tilde
S}(\mu^*K_S+B+(m-n)\mu^*H)=
h^1 O_{\tilde S}((n-m)\mu^*H)$ vanishes. One can prove this last vanishing,
using the Leray spectral sequence and our assumptions on $S$.
\qed\enddemo

We are going to apply the proposition only for smooth surfaces in $\bold P^3$
with $H$= plane
divisor, so that the assumptions on $S$ hold. In this case, we get for all $m$
an isomorphism
$$ \frac{H^0(S,\cal I_N(mH ))}{H^0(S,\cal I_C(mH ))}\to
H^0(\tilde C,\nu^*\cal O(mH )(-\tilde N)).$$

\proclaim{Theorem 3.4}  Let $S$ be a smooth surface of degree $d$ in $\bold
P^3$ and let
$C\subset S$ be a complete intersection curve
 of type $d,n$, having only $\delta$ nodes
as singularities.\par
Then $C$ is geometrically linearly
normal if and only if the nodes $N$ of $C$ impose independent conditions to the
linear system
$|(n+d-5)H|$, where $H$ is the plane divisor of $S$.\par
In particular, for $d=5$,   $C$ is geometrically linearly normal if and only if
$N$ imposes independent conditions to $|C|$, i.e., if and only if
the Severi variety ${\cal V}_{C,\delta}$ is smooth of codimension $\delta=\deg
N$.
\endproclaim
\demo{Proof}
 We use the  notation of proposition 3.3. The canonical divisor of $\tilde C$
is
$\omega=\nu^*((n+d-4)H(-\tilde N))$; on the other hand, it is clear by the
definition that
$C$ is geometrically linearly normal
if and only if $h^0(\tilde C,\nu^*(H))=4$. Now observe that on $\tilde C$,
$\nu^*(H)$ is residual to $\nu^*((n+d-5)H(-\tilde N))$. Using the previous
remark, by
Riemann-Roch one computes:
$$\multline h^0(\tilde C,\nu^*(H))=nd-p_a(nH)-1+\delta+h^0(S,\cal
I_N((n+d-5)H))- \\ -
h^0(S,\cal I_C
((n+d-5)H)).\endmultline$$
Putting $h^0(S,\cal I_N((n+d-5)H))=h^0(S,\cal O_S((n+d-5)H))-\delta+s$, with
some
computations one finds:
$$h^0(\tilde C,\nu^*(H))=4+s$$
so that $C$ is geometrically linearly normal if and only if $s=0$, i.e. if and
only if
$N$ imposes independent conditions to $|(n+d-5)H|$.\qed \enddemo

We can use the same argument of theorem 2.2 to give a partial answer to problem
(3.2).

\proclaim{Theorem 3.5}  Let $S$ be a smooth surface of degree $d\ge 5$ in
$\bold P^3$
and let $
H$ be its plane divisor; let $C\in|nH|$, $n\ge 2$, be an irreducible curve,
having only
$\delta$ nodes for singularities. If
$$ \delta< \frac{nd(n-2)}{4} $$
then $C$ is geometrically linearly normal. \endproclaim

\demo{Proof} It is very similar to the proof of theorem 2.2. We show that if
$C$ is not geometrically linearly normal, we get a contradiction.\par
Indeed, if this happens, then  the set of nodes  $N$ of $C$  does not impose
independent conditions to the curves of $|(n+d-5)H|$, by theorem
 3.4.
Take the subset $N_0\subset N$ and the rank 2 vector bundle $E$ described in
remark 2.3 and put $\delta_0=\deg( N_0)$; in this case the exact sequence is
$$0\to \cal O_S\to E\to \cal I_{N_0}((n-1)H)\to 0$$
Here $c_1(E)=\cal O_S((n-1)H)$, hence by assumption:
$$ c_1(E)^2-4c_2(E) = (n-1)^2d-4\delta_0\ge (n-1)^2d-4\delta>0 $$
so that $E$ is Bogomolov unstable; it follows that there exists a
'destabilizing' divisor $M$ for which $ h^0(S,E(-M))>0 $ and
$$(2M-c_1(E))H >0 \quad\text{i.e. }MH>({n-1\
over 2})d. \tag 7$$
Taking $M$ maximal, we may further assume that a general section of $E(-M)$
vanishes in a locus of codimension 2,whose degree $c_2(E(-M))$ must be $\ge 0;$
hence:
$$ \delta_0+M^2-(n-1)MH = c_2(E)+M^2-Mc_1(E) = c_2(E(-M)) \ge 0.$$
Now use again the assumption $\delta_0\le\delta<n(n-2)d/4$ and
observe that, by Hodge theorem, $dM^2 \le (MH)^2$; putting all
together, we arrive to the inequality:
$$ \frac{(MH)^2}{d}-(n-1)MH+{n (n-2) \over 4}d>0 \tag 8$$
from which, since   $MH<(n-2)d
/2$ yields a contradiction,  one deduces that $MH>{nd/2}$.\par
Let us now go back to the exact sequence above. $h^0(S,\cal O_S(-M))$ is $0$,
for $-MH<-(n-1)d/2< 0$ by assumptions; thus
$h^0(S,E(-M))>0$ implies $h^0(S,\cal I_{N_0}((n-1)H-M)>0$, that is, there
exists a
divisor $\Delta$ in the linear system $|(n-1)H-M|$, which contains  $\delta_0$
nodes
of the curve $C$. $\Delta$ cannot contain $C$ as a component, for
$(-H-M)H<0$. It follows, by Bezout,
$$((n-1)H-M)(nH)\ge 2\delta_0.$$
Putting all
together, one finally gets:
$$ \frac{(MH)^2}{d}-{3n-2\over 2}MH+{n(n-1)d\over 2} \ge 0. $$
Since by assumption $n\ge 2$, then we get that either  $MH\le nd/2$,
which is excluded by \thetag 8, or $MH\ge(n-1)d$; but this last inequality
yields
$((n-1)H-M)H\le 0$ while $|(n-1)H-M|$  has the
effective divisor $\Delta$ which contains $N_0\neq\emptyset$, a contradiction.
\qed
\enddemo

Using the same arrangement of proposition 2.5, one can improve the previous
statement
when $S$ is a general smooth quintic
 surface.

\proclaim{Proposition 3.6}  Let $S$ be a smooth quintic surface in $\bold P^3$,
with
Picard group $\bold Z$. Let  $C\in|nH|$ be a curve with only $\delta$ nodes as
singularities. Assume $n$ odd and
$$\delta < {5(n-1)^2 \over 4}.$$
Then $C$ is geometrically linearly normal. \endproclaim

\proclaim{Remark 3.7}\rm In the hypothesis of theorem 3.4, if $5\le n<d$ and
$C$ is also
contained in a smooth surface $S'$ of degree $n$, then one may interchange the
roles of
$n,d$ and prove that
$C$ is geometrically linearly normal in the wider range:
$$ \delta< \frac{nd(d-2)}{4}. $$\endproclaim

\bigskip

{\bf 4. THE EXAMPLES}
\medskip

Here we show the sharpness of the bounds in theorem 2.2 and theorem 3.5, for
the case
of a general smooth quintic surface $S\subset\bold P^3.$\par

{}From now on, in this section, let $S$ be a general smooth 5-ic surface of
$\bold P^3$, with Picard group $\bold Z$, generated by the plane divisor
$H$. Let $C$ be a curve in the linear system $|nH|$, $n
\ge 2$, with
$\delta$ nodes for singularities.\par
{}From theorem 2.2, we know that when
$$ \delta <  {5n (n-2 )\over 4}$$
then the Severi variety ${\cal V}_{nH,\delta}$ is smooth of
codimension $\delta$. When $n$ is odd, by proposition 2.5 the
same conclusion holds when
$$ \delta <  {5(n-1)^2 \over 4}.$$
We show with examples that these bounds are sharp. Thus we are going to
produce curves $C$ as above, with $5n(n-2)/4$ nodes for $n$ even  or
$5(n-1)^2/4$ nodes, for $n$ odd, such that the nodes do
not impose independent
conditions to the linear system $|nH|$.\par
By theorem 3.4, such a curve $C$ is a birational projection of some curve
$C'$ lying in $\bold P^4$.

\proclaim{Example 4.1} $n$ even, $n=2m$, $m\ge 3$. \par        \rm
Let $X$ be a general complete intersection surface of type $2,m$ in
$\bold P^4$ and let $X'$ be a general projection of $X$ in $\bold P^3$.
$X'$ has a double curve $Y$ of degree $m^2-m$, as one can see taking
general hyperplane sections of $X$ and $X'$.\par
Let
 $\tilde S$ be a general cone in $\bold P^4$, with vertex $V$, over our
general 5-ic surface $S$ and call $\tilde C$ the intersection of $\tilde S$
with a general
complete intersection $X$ as above. Put $C$= projection of $\tilde C$
from $V$; $C$ has degree $5n$ and it is complete intersection of $S$ and $X'$
in
$\bold P^3$, so it belongs to the linear system
$|nH|$ on $S$; moreover $C$ has nodes in the points
of $S\cap Y$, so it has a set $N$ of $\delta=5(m^2-m)=5n(n-2)/4$ nodes
and no other singularities.\par
Since $C$ is not geometrically linearly normal, it follows from theorem 3.4
that $N$ cannot
impose independent conditions to $|nH|$, so that ${\cal V}_{nH,\delta}$
is not smooth  of codimension $\delta$, in a neighbourhood of $C$.
\endproclaim

\proclaim{Proposition 4.2} The curve  $C$ constructed in the previous
example is a  singular point  of ${\cal V}_{nH,\delta}$, which is generically
smooth, of
the expected codimension $\delta$.\endproclaim
\demo{Proof} The previous construction,
 in fact, together with the proof
of theorem 3.4, shows that the tangent space of $\cal{V}_{nH,\delta}$
at $C$, that is ${H^0(S,{\cal I}_N(nH))/ (C)},$
has codimension $\delta-1$ in the tangent space of $|nH|$ at $C$:
indeed $C$ is the projection of a smooth, arithmetically normal curve
in $\bold P^4$. Hence $h^1(S,\cal I_N(nH))=1$.\par
Let $C'$ be a curve in a neighbourhood of $C$ in $\cal{V}_{nH,\delta}$ for
which the set of nodes $N'$ does not impose independent conditions to $|nH|$.
Then by semicontinuity
  $h^1(S,\cal I_{N'}(nH))=1$, so by 3.4 again, $C'$
is the projection of a curve $\tilde C'$ in $\bold P^4$ and $\tilde C'$
lives in a neighbourhood of $\tilde C$ in the Hilbert scheme of $\bold P^4$.
It follows that also $\tilde C'$ must be a smooth complete intersection of the
cone $\tilde S$ with some complete intersection surface of type $2,m$.\par
Let us compute, now, the dimension of the subvariety $\cal V$ of
$\cal{V}_{nH,\delta}$, formed by curves $C'$ at which the tangent space
of
$\cal{V}_{nH,\delta}$ has codimension $\delta-1$, in a neighbourhood
of $C$; all these curves are projection of a complete
intersection $\tilde C'$ of a cone $\tilde S$ over $S$ with a complete
intersection surface of type $2,m$ in $\bold P^4$.
If we fix the cone $\tilde S$, then  we know that such curves
$\tilde C'$ fill a variety of dimension at most
$$ h^0(\cal N_{\tilde C,\tilde S }) = h^0(\tilde C,\cal O_{\tilde
C}(2)\oplus\cal O_{\tilde C}(m))
= 14+(5m^2-10m+14)=5m^2-10m+28,   $$
so if we let
 also $\tilde S$ move, varying the vertex, then
we get a subvariety of dimension at most $5m^2-10m+32$
in the Hilbert scheme of $\bold P^4$.\par
Since  $\cal{V}_{nH,\delta}$ has dimension at least
$$  h^0(S,O_S(nH))-1-\delta= 5m^2+4  $$
and $ 5m^2+4 > (5m^2-10m+32)$ in our range, then a general element
$C"\in\cal{V}_{nH,\delta}$
does not arise from this construction. It follows that such $C"$ cannot be the
projection
of a non degenerate curve in $\bold P^4$, thus by theorem 3.4, its set of nodes
$N"$
imposes independent conditions to $|nH|$; then
the tangent space of $\cal{V}_{nH,\delta}$
at a general point has the expected codimension $\delta$, so that
$\cal{V}_{nH,\delta}$ is generically smooth of the expected codimension
$\delta$ but it is singular at the locus $\cal V$
constructed above.\qed
\enddemo

\proclaim {Remark 4.3} \rm
In the previous construction, the set of nodes $N$ of $C$ is the
intersection of $S$ with the singular locus of $X'$; one can show that,
accordingly with remar
k 2.7, the set $N$ is complete intersection, in
$\bold P^3$, of surfaces of degree $5,m,m-1$. 
\endproclaim

\proclaim{Example 4.4} $n$ odd, $n=2m+1$, $m\ge 3$. \par    \rm
We start here with a surface  $X\subset\bold P^4$ which is residue to a
plane $\pi$ with respect to a general complete
 intersection of type $2,m+1;$
$X$ is an  arithmetically Cohen-Macaulay surface.
Let $X'$ be a general projection of $X$ in $\bold P^3$.
$X'$ has a double curve $Y$ of degree $m^2$, as one can see taking
general hyperplane sections of $X$ and $X'$.\par
Let $\tilde S$ be a general cone in $\bold P^4$, with vertex $V$, over our
general 5-ic $S$ and call $\tilde C$ the intersection of $S$ with a general
surface $X$ as above. Let $C$ be the projection of $\tilde C$
from $V$; $C$ has degree $5n$ and it i
s complete intersection, in
$\bold P^3$, of $S$ and $X'$, so it belongs to the linear system
$|nH|$ on $S$; moreover $C$ has nodes in the points
of $S\cap Y$, so it has as set $N$ of $\delta=5m^2=5(n-1)^2/4$ nodes
and no other singularities.\par
Since $C$ comes from a smooth curve in $\bold P^4$, by theorem 3.4 $N$
cannot impose independent conditions to $|nH|$, so that ${\cal V}_{nH,\delta}$
is not smooth, of codimension $\delta$  in a neighbourhood of $C$.
\endproclaim

\proclaim{Proposition 4
.5} For $m\ge 5$ the curve $C$ constructed in example 4.4
is a singular point of ${\cal V}_{nH,\delta}$, which is generically smooth, of
the expected codimension $\delta$.\endproclaim
\demo{Proof} The previous construction, together with the proof
of theorem 3.4, shows, arguing as above, that the tangent space of
$\cal{V}_{nH,\delta}$
at $C$ has codimension $\delta-1$ in the tangent space of $|nH|$ at $C$.\par
 Let $C'$ be a curve in a neighbourhood of $C$ in $\cal{V}_{nH,\delta}$
for which the set o
f nodes $N'$ does not impose independent conditions to
$nH$. Then by semicontinuity  $h^1(S,\cal I_{N'}(nH))=1$, so by 3.4 again,
$C'$ is the projection of a curve $\tilde C'$ in $\tilde S$ and $\tilde C'$
lives in a neighbourhood of $\tilde C$ in the Hilbert scheme of $\tilde S$.
Let $(X\cup\pi)\cap\tilde S=\tilde C\cup\gamma$, where $\gamma$ is a plane
quintic.\par
We have an exact sequence on $\tilde C$:
$$0\to\cal N_{\tilde C,\tilde S } \to \cal O_{\tilde C}(2)\oplus\cal O_{\tilde
C}(m+1)\to
\cal
 T\to 0$$
where $\cal T$ is a torsion sheaf. Therefore, since $\tilde C$ is arthmetically
normal, we get:
$$\multline h^0\cal N_{\tilde C,\tilde S }\le h^0(\tilde C,\cal O_{\tilde
C}(2))+h^0(\tilde C,
\cal O_{\tilde C}(m+1))\le \\ \le 14+(m+1)5(2m+1)+1-p_a(\tilde C)+
h^1(\tilde C,\cal O_{\tilde C}(m+1)).\endmultline$$
Now observe that  $h^1(\tilde C,\cal O_{\tilde C}(m+1))-1$ equals the dimension
of the linear
system cut on $\tilde C$ by the quadrics through $\gamma$, therefore
$h^1(\tilde C,\cal O_
{\tilde C}(m+1))\le 9$ and since $p_a(\tilde C)=5m^2+15m+6$, we obtain
$h^0\cal N_{\tilde C,\tilde S }\le 5m^2+23$.\par
Therefore $\tilde C$ belongs to a component of the Hilbert scheme of $\tilde S$
of dimension
at most $5m^2+23$; if we let also $\tilde S$ move, varying the vertex, then
we obtain a family of curves  of dimension at most $5m^2+27$. Thus, this is an
upper bound
for the dimension of the subvariety $\cal V$ of $\cal V_{nH,\delta}$ formed by
curves
$C'$ at which the tangent space of $\cal
 V_{nH,\delta}$ has codimension $\delta-1$.
Since  $\cal{V}_{nH,\delta}$  has dimension at least
$$  h^0(S,O_S(nH))-1-\delta=10m^2+5m+3-\delta= 5m^2+5m+4 $$
and $ 5m^2+5m+4 >  5m^2+27 $ in our range, then the conclusion follows.\qed
\enddemo

\proclaim{Remark 4.6} \rm The two previous examples can be easily arranged on
general
smooth surfaces $S$ of degree $d\ge 5$ in $\bold P^3$ to provide irreducible,
non geometrically
linearly normal curves $C\subset S $ which are complete intersection of type
 $d,n$ and have
exactly
$$\delta={nd(n-2)\over 4}$$
nodes and no other singularities.\par
In particular, the bound in theorem 3.5 is sharp for all $d\ge 5$.\endproclaim

\proclaim{Proposition 4.7} Let $S$ be a general, smooth surface of degree 5 in
$\bold P^3$.
Assume:
$$ \delta\le \cases 5n(n-2)/4 & \text{  for $n$ even} \\
                    5(n-1)^2/4 & \text{ for $n$ odd} \endcases $$
Then $\cal{V}_{nH,\delta}$ is non empty and it has at least one generically
smooth component of codimen
sion $\delta$ in $\bold P(|nH|)$. \endproclaim
\demo{Proof}  For $\delta= n(n-2)/4$,  $n$ even, or $\delta=(n-1)^2/4$, $n$
odd, the statement
follows by proposition 4.2 and proposition 4.4.\par
For smaller $\delta$, one can argue as in the proof of theorem 1.1. For $n$
even, start
with a general curve  $C\in\cal{V}_{nH,n(n-2)/4}$; by proposition 4.2 or
proposition 4.5,
the nodes of $C$ impose independent
conditions to the curves of $|nH|$, so we may smooth them independently, one by
one, getting
at
any step curves with only nodes for singularities.
The conclusion follows from theorem 2.2.\par
A similar argument works for $n$ odd.
\qed \enddemo

 \end